\documentclass{article}

\PassOptionsToPackage{numbers, compress}{natbib}
\usepackage{cite}
\usepackage[pagebackref]{hyperref}

\usepackage[preprint]{neurips_2024}
\usepackage[utf8]{inputenc}
\usepackage[T1]{fontenc}
\usepackage{hyperref}
\usepackage{url}
\usepackage{booktabs}
\usepackage{amsfonts}
\usepackage{nicefrac}
\usepackage{microtype}
\usepackage{xcolor}

\usepackage{csquotes}
\usepackage{CJKutf8}
\usepackage{graphicx}
\usepackage[most]{tcolorbox}
\usepackage{cleveref}
\usepackage{multirow}

\usepackage{array}
\usepackage{booktabs}

\usepackage{amsmath}
\usepackage{amssymb}
\newtheorem{theorem}{Theorem}[section] 
\newtheorem{proposition}[theorem]{Proposition} 
\newtheorem{lemma}[theorem]{Lemma}

\title{The Hurwitz existence problem and \\the  prime-degree conjecture: \\
A computational perspective}

\author{
    Yiru Wang$^{1}$, Bingqian Li$^{2}$, Yi Zhou$^{2}$, Zhiqiang Wei$^{3}$, Yu Ye$^{1}$, Yiqian Shi$^{1}$, Bin Xu$^{1}$ \\
    $^1$School of Mathematical Sciences, University of Science and Technology of China(USTC)\\
    $^2$National Engineering Laboratory for Brain-inspired Intelligence Technology and Applications,\\School of Information Science and Technology,  USTC\\
    $^3$School of Mathematics and Statistics, Henan University
}

\newcommand\RS{\mathbb{CP}^1}

\begin{document}

\maketitle

\begin{abstract}
We investigate the Hurwitz existence problem from a computational viewpoint. Leveraging the symmetric-group algorithm by Zheng and building upon implementations originally developed by Baroni, we achieve a complete and non-redundant enumeration of all non-realizable partition triples for positive integers up to $31$. These results are further categorized into four types according to their underlying mathematical structure; it is observed that nearly nine-tenths of them can be explained by known theoretical results. As an application, we verify the prime-degree conjecture for all primes less than $32$. In light of the exponential memory growth inherent in existing computational approaches---which limits their feasibility at higher degrees---we propose a novel software architecture designed to stabilize memory usage, thereby facilitating further detection of exceptional cases in the Hurwitz existence problem. The complete dataset of non-realizable partition triples, along with our implementation, will been made public on GitHub.
\end{abstract}

\noindent MSC 2010: 57M12 (primary), 20B35, 05A17.
\section{Introduction}

The classical Hurwitz existence problem concerns the realizability of prescribed branch data by branched coverings of the Riemann sphere $\RS$ \cite{Hurwitz:1891, EKS1984}. Given a collection $\Lambda$ of $k \geq 3$ nontrivial partitions of a positive integer $d$, one asks whether there exist a compact connected Riemann surface $S$ together with a branched covering $f: S \rightarrow \RS$ whose branch data agree with $\Lambda$. When such a map $f$ exists, we say that $\Lambda$ is {\it realizable}. Any realizable collection must satisfy an algebraic necessary condition, given by the Riemann-Hurwitz formula \cite[\S 7.2.1, p.~101]{Donaldson2011}. We refer to this as the {\it Riemann-Hurwitz condition}. A collection satisfying this condition but not realizable will be called {\it exceptional}.

Hurwitz \cite{Hurwitz:1891} reduced this problem to an equivalent one of constructing $k$ suitable permutations in the symmetric group $S_d$ associated with the collection $\Lambda$. See \cite[\S\,4.2]{Donaldson2011} for a proof of this reduction. Along this line, Mednykh \cite{Mednykh1984} gave a conceptual solution to the Hurwitz existence problem in the sense that he obtained a complicated formula for the number of equivalence classes of branched covers realizing $\Lambda$. However, it seems difficult in practice to determine the realizability of a given collection $\Lambda$ from Mednykh's formula. Therefore, beyond the Riemann-Hurwitz condition, the problem remains highly intricate and is known to encode deep interactions between complex analysis, algebraic geometry, and the representation theory of symmetric groups. The reader may consult \cite{SX2020, Pet2020, BP2024} and the references therein for up-to-date theoretical results on the Hurwitz existence problem.

A particularly appealing case is the prime-degree conjecture proposed by Edmonds-Kulkarni-Stong \cite{EKS1984}, which asserts that whenever the degree $d$ is prime, any collection $\Lambda$ of $k$ nontrivial partitions of $d$ satisfying the Riemann-Hurwitz condition must be realizable. They also reduced the conjecture to the case $k=3$.
This conjecture was given strong support by the work of Zheng \cite{Zheng2006}, who reformulated the realizability problem in terms of generating functions involving secondary partitions and characters of $S_d$, and confirmed the conjecture for primes less than $20$. Building on Zheng's algorithm, Baroni \cite{Baroni2024, BP2024} carried out further verifications, confirming the conjecture for all primes $d \le 29$. However, Baroni's computational approach suffers from a fundamental obstacle: it exhibits exponential memory growth, which limits the applicability for larger degrees.

In this work we revisit the Hurwitz existence problem from a computational perspective. Building on Zheng's algorithm, we undertake a new computational architecture of the Hurwitz existence problem that overcomes this barrier. For simplicity, we restrict ourselves to the case where $\Lambda$ consists of exactly three partitions.

Our contributions in this manuscript are threefold as follows.

First, by extending Baroni's implementation and conducting computations on high-performance clusters, we obtain a complete and non-redundant enumeration of all non-realizable partition triples for every $d \le 31$, both prime and composite. This improves the known range from $d \le 29$ and provides for the first time a full catalog of all exceptional data up to degree $31$. As an immediate consequence, our enumeration verifies the prime-degree conjecture for all primes less than 32.

Second, the resulting data allow us to perform a systematic structural classification of the exceptional cases. We find that approximately 90\% of all non-realizable triples fall into two families of the four---each corresponding to explicit theoretic constraints. Thus, the overwhelming majority of exceptional triples admit a clear conceptual explanation rather than arising from sporadic computational accidents.

Third, we introduce a new software architecture, which is inspired by Baroni's  work, designed to overcome the memory limitations inherent in previous approaches. Our design employs a two-layer architecture: a Python layer manages workflow and I/O, while computationally intensive kernels are JIT-compiled for performance. This separation, combined with hash-table data structure and batch processing, reduces peak memory consumption dramatically, enabling computations that were previously unattainable. This reduces the peak memory requirement for $d = 32$ from multi-terabyte scales to under 400 GB in practice. Our design also emphasizes sustainability and resumability by providing detailed logging and checkpointing mechanisms suitable for long-running computations.

The complete dataset of exceptional partition triples for $d \le 31$ provides a refined picture of the exceptional landscape, and we hope that both the full dataset and its structural organization will facilitate future theoretical work on the Hurwitz existence problem.

The paper is organized as follows. Section 2 recalls Zheng's character-theoretic algorithm and establishes notation. Section 3 describes the false-positive removal method. Section 4 classifies these triples into four families. Section 5 describes the new computational architecture. All exceptional data and code will be deposited in the repository. 

\section{Preliminaries}\label{sec2}

Throughout this section we recall the algorithm for realizability introduced by Zheng \cite{Zheng2006}. We fix a positive integer $d$, and let $\mathcal{P}(d)$ denote the set of all partitions of $d$.
For a partition $\lambda=(\lambda_{1},\dots,\lambda_{r})\in\mathcal{P}(d)$, we use the following notation:
\begin{itemize}
\item $S^{\lambda}$ denotes the irreducible representation of the symmetric group $S_{d}$ corresponding to $\lambda$;
\item $\chi^{\lambda}$ denotes the irreducible character of $S_{d}$ afforded by $S^{\lambda}$;
\item $C_{\lambda}\subseteq S_{d}$ denotes the conjugacy class consisting of permutations whose cycle-type is $\lambda$.
\end{itemize}

\subsection{Secondary partitions}
Following Zheng, a secondary partition of $d$: $\omega=[\mu_{1},\dots,\mu_{n}]\vdash\vdash d,$
where $\mu_{j}\vdash d_{j}\ ,\ d_j\in\mathbb{Z}_{>0}\ ,\forall\ j$ such that $d_{1}+\dots+d_{n}=d$. For such $\omega$, let $n_{1},\dots,n_{m}$ be the multiplicities of the distinct components among the set of $\{\mu_{i_j}\, ,\ j=1,\dots, m\}$, so that $n_{1}+\dots+n_{m}=n$.

\subsection{The coefficients $r(\omega)$}
Define
\begin{equation}
r(\omega)\;=\;\frac{(-1)^{n-1}}{n}\,
\frac{n!}{n_{1}!\dots n_{m}!}\,
\prod_{j=1}^{n}\left(\frac{\dim S^{\mu_{j}}}{d_{j}!}\right)^2\;.
\label{eq:romega}
\end{equation}

\subsection{The polynomials $s(\omega\,;\boldsymbol{t})$}
Let $k\geq 1$ be the number of branch points.
For $D = (\mu_{1}, \ldots, \mu_{k})$, denote \[\mathbf{t}^{D} = \mathbf{t}_{1}^{\mu_{1}} \cdots \mathbf{t}_{k}^{\mu_{k}}.\]
For $\mu = [a_{1}, \ldots, a_{r}]$, $i\in\{1,\dots,k\}$, denote \[\mathbf{t}_{i}^{\mu} = t_{i,a_{1}} \cdots t_{i,a_{r}}.\]
Given a secondary partition $\omega=[\mu_{1},\dots,\mu_{n}]$ with $\mu_{j}\vdash d_{j}$, define
\begin{equation}
s(\omega\,;t_{i1},\dots t_{id})\;=\;
\prod_{j=1}^{n}
\left(
\sum_{\nu\ \vdash d_{j}}\frac{\chi^{\mu_{j}}(C_\nu)\,|C_\nu|}{\dim S^{\mu_{j}}}
\textbf{t}_{i}^{\nu}\;\right).
\label{eq:somega}
\end{equation}

\subsection{Realizability polynomial}
Define
\begin{equation}
F_{d,k}(\boldsymbol{t})\;=\;
\sum_{\omega\vdash\vdash d}\,
r(\omega)\prod_{i=1}^{k}s(\omega\,;t_{i1},\dots,t_{id}).
\label{eq:Fdk}
\end{equation}
Zheng proved that $F_{d,k}$ contains all information needed to determine realizability of branch data over $k$ branch points.

\begin{theorem} {\rm (Zheng \cite[Theorem 4]{Zheng2006})}
Let $\Lambda=(\lambda_{1},\dots,\lambda_{k})$ be a $k$-tuple of nontrivial partitions of $d$.
Then $\Lambda$ is realizable by a branched covering of compact Riemann surfaces if and only if the coefficient of the monomial
\[
\textbf{t}_{1}^{\lambda_{1}}\cdots \textbf{t}_{k}^{\lambda_{k}}
\]
in $F_{d,k}(\boldsymbol{t})$ is nonzero.
\end{theorem}

Thus the Hurwitz problem reduces to detecting which coefficients in the polynomial $F_{d,k}$ vanish.

\section{Eliminating False Positives}

The detection of non-realizable branch data in Zheng's symmetric-group framework relies on evaluating the polynomial
\[
F_{d,k}(\textbf{t})\;=\;
\sum_{\omega\vdash\vdash d}\,r(\omega)\prod_{i=1}^{k}s(\omega;t_{i1},\dots,t_{id})
\]
where $k=3$ for three branch data case. That means, $\forall \omega\vdash\vdash d$, we should estimate 
\begin{align}
    &\frac{(n-1)!}{n_{1}!\dots n_{m}!}\,\cdot
    \left(\prod_{j=1}^{n}\frac{(\dim S^{\mu_{j}})^2}{(d_{j}!)^2}\right)\,\cdot
    \prod_{i=1}^{3}\prod_{j=1}^{n}
    \left(
    \sum_{\nu\vdash d_{j}}\frac{\chi^{\mu_{j}}(C_\nu)\,|C_\nu|}{\dim S^{\mu_{j}}}
    \right)
    \notag\\
    =&\frac{(n-1)!}{n_{1}!\dots n_{m}!}\,\cdot
    \prod_{j=1}^{n}\left(\frac{(\dim S^{\mu_{j}})^2}{(d_{j}!)^2}\,\cdot
    \frac{\left(\sum_{\nu\vdash d_{j}}\chi^{\mu_{j}}(C_\nu)\,|C_\nu|\right)^3}{(\dim S^{\mu_{j}})^3}\right)
    \notag\\
    =&\frac{\frac{n!}{n_{1}!\dotsm n_{m}!}\cdot \left(\sum_{\nu\vdash d_{j}}\chi^{\mu_{j}}(C_\nu)\,|C_\nu|\right)^3}{n \cdot \prod_{j=1}^{n}\dim S^{\mu_{j}}\cdot (d_j!)^2}
    \label{eq:long}
\end{align}
$\Lambda$ is exceptional precisely when the coefficient of its corresponding monomial in $F_{d,3}$ vanishes.

In actual computation, however, to both avoid the explosion of intermediate coefficients and guarantee exact arithmetic without floating-point errors, our implementation evaluates these coefficients modulo a fixed prime $p$. Thus a genuine zero coefficient certainly appears as zero mod $p$, but the converse need not hold: a nonzero integer coefficient may also vanish modulo $p$. Such spurious outputs are false positives, and removing them is essential for producing a complete and non-redundant list of exceptional triples.

To eliminate false positives rigorously, we employ the following procedure, justified by a bound on the numerator of coefficients.

\subsection{Bounding the coefficients of $F_{d,3}$}

For a secondary partition $\omega=[\mu_{1},\dots,\mu_{n}]\vdash\vdash d$, the polynomial
\[
s(\omega;\boldsymbol{t})\;=\;
\prod_{j=1}^{n}
\left(
\sum_{\nu\ \vdash d_{j}}\frac{\chi^{\mu_{j}}(C_\nu)\,|C_\nu|}{\dim S^{\mu_{j}}}
\textbf{t}_{i}^{\nu}\;\right).
\label{eq:somega}
\]
is a product of $n$ polynomials.

First, we shall estimate the numerator part of the above unexpanded form.

Using the first-orthogonality relation
\[
\sum_{g\in S_{d_{j}}}|\chi^{\mu_{j}}(g)|^{2}=d_{j}!,
\qquad
\sum_{g\in S_{d_{j}}}1=d_{j}!,
\]
and the Cauchy-Schwarz inequality, we have
\[
\sum_{g\in S_{d_{j}}}|\chi^{\mu_{j}}(g)|
\le
\Bigl(\sum_{g\in S_{d_{j}}}|\chi^{\mu_{j}}(g)|^{2}\Bigr)^{1/2}
\Bigl(\sum_{g\in S_{d_{j}}}1\Bigr)^{1/2}
=d_{j}!.
\]
Since $|\chi^{\mu_{j}}(g)|$ is constant on conjugacy classes, we write
$\chi^{\mu_{j}}(C_{\nu})$ for $g\in C_{\nu}$ ($\nu\vdash d_{j}$) and obtain an upper bound
\[
\sum_{\nu\vdash d_{j}}\,|\chi^{\mu_{j}}(C_{\nu})||C_{\nu}|
\le
d_{j}!.
\]

Next, we derive an upper bound for the numerators of the rational coefficients after expressed with a common denominator.

Fix a secondary partition $\omega=[\mu_{1},\dots,\mu_{n}]\vdash\vdash d$. When $k=3$, the un-estimated part of $|r(\omega)\prod_{i=1}^{3}s(\omega;t^{(i)})|$ can be written as
\begin{equation}
\frac{\frac{n!}{n_{1}!\dotsm n_{m}!}}{n \cdot\prod_{j=1}^{n}\dim S^{\mu_{j}}\cdot (d_j!)^2}.
\label{eq:bound}
\end{equation}
Using $n_1+\dots +n_m=n$, we have that \[\frac{n!}{n_{1}!\dotsm n_{m}!}\in \mathbb{Z}_{>0}, \]
and \[\frac{n!}{n_{1}!\dotsm n_{m}!}\le n!<d\,!.\]
Since \[\dim S^{\mu_{j}} \le \ d_j!\ ,\] together with \[d_1+\dots +d_n = d\ ,\] it follows that the numerator of \eqref{eq:long} can be bounded by $(d!)^4$, while the denominator can be bounded by $n \cdot (d!)^3$.

Now, we consider the sum over all $\omega$. In order to apply the Chinese Remainder Theorem(CRT), we must first find a common denominator.

Because \[\dim S^{\mu_{j}} | \ d_j!\ ,\] we may take $(d!)^4$ as a common denominator. With this denominator, the numerator of each summand corresponding to a $\omega$ is bounded by $(d!)^8$. Therefore, if there are $p_{2}(d)$ secondary partitions $\omega\vdash\vdash d$, summing over the number of $\omega$ we obtain
\begin{lemma}[Numerator bound]
\label{lem:num-bound}
After using a common denominator as $(d!)^4$, every coefficient of the polynomial $F_{d,3}$ is a rational number whose numerator (in absolute value) is at most
\[
p_{2}(d)\,\cdot (d!)^8.
\]
\end{lemma}

\subsection{Testing modulo sufficiently many primes}

Let
\[
M_{d}\;:=\;p_{2}(d)\cdot(d!)^{8}.
\]
If we choose primes $p_{1},\dots,p_{N}$ whose product exceeds $M_d$, then the following holds:

\begin{proposition}[Chinese remainder theorem for false positives]
Suppose a set of three branch data $\Lambda$ gives zero coefficient in $F_{d,3}$ modulo each of the primes $p_{1},\dots,p_{N}\ (\ p_i \ne p_j\;,\  \forall i \ne j\ )$.
If
\[
p_{1}p_{2}\cdots p_{N}\;>\;M_{d},
\]
then the true integer coefficient of $\Lambda$ in $F_{d,3}$ is zero.  
In particular, $\Lambda$ is genuinely non-realizable.
\end{proposition}

Thus by running the exceptional-checking program under a list of primes whose product exceeds $M_{d}$, any triple that survives all tests is confirmed to be a true exceptional triple.

\subsection{Summary on False Positives for $d\le 30$}

We initially compute all the 'exceptional data' under the prime number 1,000,000,007. Then, using the aforementioned verification method, we identify all the false positives for $d\le 30$, which are not genuine exceptional data, as shown in the table below.

\begin{tabular}{l p{3.9cm} p{3.2cm} p{2.8cm}}
\toprule
\multicolumn{1}{c}{$d$} &\multicolumn{1}{c}{$\lambda_1$} &\multicolumn{1}{c}{$\lambda_2$} &\multicolumn{1}{c}{$\lambda_3$} \\
\midrule
25    & [ 8 8 2 2 2 1 1 1 ] &[ 7 4 4 2 2 2 2 1 1 ]& [ 11 5 3 2 2 2 ]  \\
\midrule
26    & [ 10 5 3 2 1 1 1 1 1 1 ] &[ 7 6 2 2 2 2 2 1 1 1 ]& [ 10 9 4 3 ]\\
\midrule
28&[ 4 3 3 3 3 3 2 2 2 2 1 ]& [ 15 5 4 2 1 1 ] &[ 11 7 4 3 3 ]\\
&[ 10 4 4 2 2 1 1 1 1 1 1 ]&[ 7 5 3 3 3 3 3 1 ] & [ 11 10 3 1 1 1 1 ]\\
&[ 7 6 3 2 2 2 1 1 1 1 1 1 ]& [ 4 3 3 3 3 2 2 2 2 2 2 ] &[ 19 4 3 1 1 ]\\
&[ 7 5 3 2 2 2 1 1 1 1 1 1 1 ] &[ 6 6 4 3 3 2 2 1 1 ]&[ 18 4 4 2 ]\\
\midrule
30&[ 12 8 4 1 1 1 1 1 1 ]&[ 6 6 5 5 5 3 ] & [ 19 7 4 ]\\
&[ 11 6 6 1 1 1 1 1 1 1 ] &[ 11 5 4 4 2 2 2 ] &[ 7 6 4 4 4 3 2 ]\\
&[ 8 4 4 4 2 2 2 1 1 1 1 ]& [ 14 5 3 2 2 2 1 1 ]& [ 14 6 6 1 1 1 1 ]\\
& [ 4 4 4 4 3 3 2 2 1 1 1 1 ]& [ 5 4 4 4 3 3 2 2 2 1 ]& [ 14 10 1 1 1 1 1 1 ]\\
& [ 11 4 3 2 2 2 1 1 1 1 1 1 ]& [ 9 7 6 3 3 1 1 ]&[ 14 5 4 2 2 2 1 ]\\
&[ 6 6 3 3 2 2 1 1 1 1 1 1 1 1 ] &[ 16 5 4 2 1 1 1 ]&[ 14 8 4 3 1 ]\\
\bottomrule
\end{tabular}

\section{Classification of Exceptional Triples}
\label{sec:classification}

Having obtained the complete and cleaned catalogue of exceptional partition triples for degrees $d\le 31$, we organise them into four disjoint types.

\paragraph{Type 0.}
\hfill

The sum of the lengths of the three partitions is less than $d+2$.
Note that if a collection of three non-trivial partitions of $d$ is realizable 
and the sum of the lengths of these three partitions is less than $d+2$,
then this collection is realizable by a branched cover from a compact Riemann surface with positive genus
to $\RS$.

\paragraph{Type I.}
\hfill

The sum of the lengths equals $d+2$, and exists one partition which contains a unique part greater than $1$. Both G. Boccara \cite{Boccara1982} and Song-Xu \cite{SX2020} proved more general theorems, predicting that all the data of Type I are exceptional.

\paragraph{Type II.}
\hfill

If the triple is not of Type 0 or I, and there exists a common divisor $c>1$ of all the parts of two partitions, while the third partition cannot be combined into $c$ partitions $\mu_1, \dots ,\mu_c$, such that \[\mu_i \vdash \frac{d}{c}, \forall i=1,\dots c,\] then the triple is assigned to \emph{Type II}.  The main theorem in Wei-Wu-Xu \cite{WWX2024} showed such data are exceptional.

\paragraph{Type III.}
\hfill

All exceptional triples that do not satisfy the previous conditions are collected in Type III.

\hfill

We tag every exceptional triple with one of the four labels Types 0, I, II, III. Empirically, Types I and II account for roughly 90\% of all exceptions. In fact, the conditions defining Types I and II are known to be sufficient to guarantee exceptionality, providing theoretical justification for these two classes. The authors have no theoretical understanding about the data of Type 0 or III. 

\section{Computational Architecture}

Two aspects of formula~\eqref{eq:Fdk} have a critical impact on performance:
\begin{enumerate}
    \item The number of secondary partitions grows extremely rapidly with $d$, reaching several million already at $d=32$.
    \item Expanding each $s(\omega;\mathbf{t})$ generates a large number of monomials. Each coefficient requires sums over the partitions, making the computation time-consuming, and the sheer number of coefficients makes storage memory‑intensive.
\end{enumerate}
Our computational architecture mitigates these bottlenecks through batched processing, a hash‑table‑based data structure, and JIT‑compiled parallel kernels for compute‑intensive loops.

\subsection{Compatibility requirements}\label{compatible}

Let $\mathcal{P}(d)$ denote the set of all partitions of $d$, ordered in reverse lexicographic order. For each $\lambda\in\mathcal{P}(d)$, we define its length
\[\ell(\lambda)=\text{the number of parts of }\lambda.\]

A triple $(\lambda_1,\lambda_2,\lambda_3)$ with $\lambda_i\vdash d$ that is realized as the branch data of a branched covering must satisfy the numerical constraints imposed by the Riemann-Hurwitz formula.  For a branch point of cycle type~$\lambda$, the ramification contribution equals $d - \ell(\lambda)$.  Substituting these contributions into the Riemann-Hurwitz formula yields
\begin{equation}\label{RH}
    2g(S)-2 = -2d + \sum_{i=1}^{3}\bigl(d-\ell(\lambda_i)\bigr).
\end{equation}
Imposing the necessary condition $g(S)\ge 0$ gives
\begin{equation}\label{ineq}
    \ell(\lambda_1)+\ell(\lambda_2)+\ell(\lambda_3)\le d+2.
\end{equation}
Taking~\eqref{RH} modulo~$2$, we also obtain the parity condition
\begin{equation}\label{pmod}
    \ell(\lambda_1)+\ell(\lambda_2)+\ell(\lambda_3)\equiv d \pmod{2}.
\end{equation}

\eqref{ineq} and \eqref{pmod} are the constraints compatitable with the Riemann-Hurwitz formula.

\hfill

The subsections below follow the execution flow of the algorithm. Each subsection corresponds to a major stage in the computational pipeline.

\subsection{Precomputation}

Before the coefficient computation begins, we construct several combinatorial and representation–theoretic tables associated with the symmetric groups 
$S_1,\dots, S_d$.

\begin{itemize}
    \item Partition table.
    \hfill
    
    For each $0\le n\le d$, we enumerate all partitions of $n$ in reverse lexicographic order, together with their lengths, offsets, and cumulative counts. This provides constant-time access to every partition and enables binary search within each layer for $n$.
    \item Conjugacy-class sizes.
    \hfill
    
    For each partition \(\nu\vdash n\) we compute the size of the conjugacy class \(C_\nu\).
    \item Partition multiplication table.
    \hfill
    
    We construct a table encoding the product of monomials. Concretely, given two monomials corresponding to partitions $\lambda$ and $\mu$, their product corresponds to the sorted multiset of $\lambda\cup\mu$. The table stores the offset of this resulting partition, enabling constant-time lookup during polynomial multiplication.
    
    \item Irreducible character table.
    \hfill
    
    The values \(\chi^\lambda(C_\nu)\) for all \(\lambda,\nu\vdash n\), \(0\le n\le d\), are computed recursively using the Murnaghan–Nakayama rule via rim–hook removal, implemented with all arithmetic performed modulo a fixed prime \(p\).
    
    \item Factorials and inverse factorials modulo a prime $p$.
    \hfill
    
    The values $n!$ (mod $p$) and $(n!)^{-1}$ (mod $p$) for \(0\le n\le d\) are computed.
\end{itemize}

For $d=32$, the full precomputation of these tables takes only 7 seconds in our experiment.

\subsection{Batch Computation for $r(\omega)$ and $s(\omega)$}

To compute $r(\omega)$ and $s(\omega;\textbf{t})$ efficiently, we process the primary partitions in batches, where each batch contains $B \in \mathbb{N}$ partitions (typically $100 \le B \le 200$). 

For each primary partition $\lambda$, all secondary partitions $\omega$ are enumerated recursively by subdividing each primary partition; in practice, to control stack depth, we use stack to replace recursion. For every $\omega$, we compute $r(\omega)$ and a compressed form of the $s(\omega)$, without storing $\omega$. To avoid memory blow-up, polynomials $s(\omega)$ are never accumulated globally; instead, each batch is streamed directly to a disk-backed storage. 

A batch size $B=150$ for $d=32$ occupies approximately $40$GB on disk, allowing computations for $d\ge 30$ without terabyte-scale memory.

\subsection{Accumulation for Each Fixed $\lambda_i$}

To test the realizability of triples with first component $\lambda_i$, we stream each batch from disk and update partial sums according to \eqref{eq:Fdk}.

Instead of storing all $(j,k)$ coefficients explicitly, we adopt a cache-friendly segmented hash table, with each segment roughly corresponding to a CPU cache line. The hash table uses the Cantor pairing function to map each key $(j,k)$ to a 64-bit integer. Collisions are resolved via quadratic probing, which searches for an available slot, backed by a small overflow buffer for robustness. 

Each $\lambda_i$ owns an independent accumulator, enabling embarrassingly parallel work distribution and high memory efficiency.

\subsection{Zero-Detection and Exceptional Triples}

Within the hash-table structure, a triple is classified as exceptional if:
\begin{enumerate}
    \item it satisfies the necessary conditions mentioned in ~\ref{compatible}~, and
    \item its key is either absent from the hash table or has value zero.
\end{enumerate}

\section*{Acknowledgements} 
B.X. would like to express his sincere gratitude to Yufei Bai, Xiaoyang Chen, Tielin Dai, Dun Liang, Sicheng Lu, Xin Nie, Zijin Peng, Qingyu Xu, Junyi Yang, Wenyan Yang, Xiaoya Zhai and Youliang Zhong for many valuable and stimulating conversations during the course of this project. The authors are also grateful to Chao Wu at the Supercomputing Center of USTC for kindly answering their questions on the supercomputing system there.

Both B.L. and Y.Z. are supported in part by Anhui Province Science and Technology Tackle Plan 
Project (No. 202423k09020008). B.X. is supported in part by the Project of Stable Support for Youth Team in Basic Research Field, CAS (Grant No. YSBR-001) and NSFC (Grant No. 12271495). The numerical calculations in this manuscript have been carried out on the supercomputing system in the Supercomputing Center of USTC.

\clearpage
{\small
\bibliographystyle{unsrtnat}
\bibliography{references}
}

{
\raggedright
YIYU WANG, YE YU, YIQIAN SHI, BIN XU \\
SCHOOL OF MATHEMATICAL SCIENCES, USTC \\
HEFEI 230026 CHINA\\
\texttt{wyr221072@mail.ustc.edu.cn\\  yeyu@ustc.edu.cn\\  yqshi@ustc.edu.cn\\ bxu@ustc.edu.cn}

\vspace{1cm}

BINGQIAN LI, YI ZHOU \\
NATIONAL ENGINEERING LABORATORY FOR BRAIN-INSPIRED INTELLIGENCE TECHNOLOGY
AND APPLICATIONS\\
SCHOOL OF INFORMATION SCIENCE AND TECHNOLOGY,  USTC\\
HEFEI 230026 CHINA\\
\texttt{bqli315@gmail.com\\ yi\_zhou@ustc.edu.cn}

\vspace{1cm}

ZHIQIANG WEI\\
SCHOOL OF MATHEMATICS AND STATISTICS\\
HENAN UNIVERSITY\\
KAIFENG 475004 CHINA\\
\texttt{weizhiqiang15@mails.ucas.edu.cn}
}

\end{document}